\documentclass[12pt]{amsart}
\usepackage{amsmath}
\usepackage{amssymb}
\usepackage{bm}
\usepackage{graphicx}
\usepackage{verbatim}
\usepackage{enumitem}
\setlist[enumerate,1]{leftmargin=*}
\setlist[itemize,1]{leftmargin=*}
\usepackage{mathrsfs}
\usepackage[dvipsnames]{color}
\usepackage{todonotes}
\usepackage[pdftex,pdfstartview=FitH,pdfborderstyle={/S/B/W 1},%
colorlinks=true, linkcolor=blue, urlcolor=blue, citecolor=blue,%
pagebackref=true]{hyperref}
\usepackage[nobysame,alphabetic,initials,msc-links]{amsrefs}
\usepackage{imakeidx}

\newtheorem{Lemma}{Lemma}[section]
\newtheorem{Theorem}[Lemma]{Theorem}
\newtheorem{Proposition}[Lemma]{Proposition}
\newtheorem{Corollary}[Lemma]{Corollary}

\theoremstyle{definition}

\newtheorem{Definition}[Lemma]{Definition}

\numberwithin{equation}{section}

\newcommand {\N} {{\mathbb N}}

\newcommand {\R} {{\mathbb R}}

\newcommand {\Z} {{\mathbb Z}}
\newcommand {\absolute}[1] {\left| {#1} \right|}
\newcommand {\norm}[1] {\left\| {#1} \right\|}

\DeclareMathOperator {\SL} {SL}
\DeclareMathOperator {\SO} {SO}

\DeclareMathOperator {\Kak} {Kak}

\DeclareMathOperator {\ad} {ad}

\DeclareMathOperator {\bu} {\mathbf{u}}
\DeclareMathOperator {\ba} {\mathbf{a}}
\DeclareMathOperator {\bx} {\mathbf{x}}
\DeclareMathOperator {\bw} {\mathbf{w}}
\DeclareMathOperator {\by} {\mathbf{y}}

\newcounter{c}

\newcounter{C}
\newcommand{\Cnew}{%
\refstepcounter{C}%
\ensuremath{C_{\theC}}}
\newcommand{\Cold}[1]{\ensuremath{C_{\ref{#1}}}}

\newcounter{e}
\newcommand{\enew}{%
\refstepcounter{e}%
\ensuremath{\epsilon_{\thee}}}
\newcommand{\eold}[1]{\ensuremath{\epsilon_{\ref{#1}}}}

\newcounter{de}
\newcommand{\denew}{%
\refstepcounter{de}%
\ensuremath{\delta_{\thede}}}
\newcommand{\deold}[1]{\ensuremath{\delta_{\ref{#1}}}}

\newcounter{d}

\newcounter{R}
\newcommand{\Rnew}{%
\refstepcounter{R}%
\ensuremath{R_{\theR}}}
\newcommand{\Rold}[1]{\ensuremath{R_{\ref{#1}}}}

\newcounter{r}

\newcounter{K}
\newcommand{\Knew}{%
\refstepcounter{K}%
\ensuremath{K_{\theK}}}
\newcommand{\Kold}[1]{\ensuremath{K_{\ref{#1}}}}

\newcounter{k}

\newcounter{w}

\newcounter{ps}

\newcounter{th}

\newcounter{n}

\makeindex[title=Index of sets and constants,columns=2]

\begin{document}
\title{Time Change for unipotent Flows and Rigidity}

\author{Elon Lindenstrauss}
\thanks{E.L. acknowledges support by ERC 2020 grant HomDyn (grant no. 833423)}
\address{E.L.: Einstein Institute of Mathematics, Edmond J. Safra Campus,
The Hebrew University of Jerusalem,
Givat Ram, Jerusalem, 9190401, Israel}
\email{elon@math.huji.ac.il }
\author{Daren Wei}
\thanks{D.W. acknowledges support by ERC 2020 grant HomDyn (grant no. 833423)}
\address{D.W.: Einstein Institute of Mathematics, Edmond J. Safra Campus,
The Hebrew University of Jerusalem,
Givat Ram, Jerusalem, 9190401, Israel}
\email{Daren.Wei@mail.huji.ac.il}
\date{\today}

\maketitle
\begin{abstract}
We prove a dichotomy regarding the behavior of one-parameter unipotent flows on quotients of semisimple lie groups under time change. We show that if $u^{(1)}_t$ acting on $G_{1}/\Gamma_1$ is such a flow it satisfies exactly one of the following:
\begin{enumerate}
    \item The flow is loosely Kronecker, and hence isomorphic after an appropriate time change to any other loosely Kronecker system.
    \item The flow exhibits the following rigid behavior: if the one-parameter unipotent flow $u^{(1)} _ t$ on $G_1/\Gamma_1$ is isomorphic after time change to another  such flow $u^{(2)} _ t$ on $G_2/\Gamma _ 2$, then $G_1/\Gamma_1 $ is isomorphic to $G_2/ \Gamma_2$ with the isomorphism taking $u^{(1)} _ t$ to $u^{(2)} _ t$ and moreover the time change is cohomologous to a trivial one. 
\end{enumerate}
The full details will appear in \cite{Lindenstrauss2022}.
\end{abstract}

\section{Introduction}

Unipotent flows have some striking rigidity properties: for instance, every orbit is recurrent \cite{margulis1971}, and every orbit closure as well as every invariant measure is homogeneous \cites{ratner1990measure,ratner1991measure,ratner1991topological}. 

In particular, they are known to be rigid with respect to measurable isomorphisms.
To fix notations, for $i=1,2$, let $G_i$ be the identity component of a real semisimple linear algebraic group without compact factors, and let $\Gamma _ i < G_i$ be a lattice. Let $\mathcal B_i$ be the Borel $\sigma$-algebra on $G _ i / \Gamma _ i$, let $m_i$ on $G _ i / \Gamma _ i$ be the normalized Haar measure, and suppose that $u^{(i)}_t$ are one parameter unipotent subgroups of $G _ i$ respectively. This means that $u^{(i)}_t = \exp (t \mathbf u^{(i)})$ with $\mathbf u^{(i)}$ a nilpotent element of the Lie algebra $\mathfrak g_i$ of $G_i$ (considered as a matrix). 
In \cite{ratner1982rigidity}, Ratner showed that if $G_1= G_2 = \SL_2(\R)$ and $\psi: G _ 1 / \Gamma _ 1 \to G _ 2 / \Gamma _ 2$ is a measurable isomorphism between the flows $(G _ 1 / \Gamma _ 1, \mathcal B_1, m_1, u^{(1)}_t)$ and $(G _ 2 / \Gamma _ 2, \mathcal B_2, m_2, u^{(2)}_t)$, then $\psi$ is automatically given by a very simple (and algebraic) form. Concretely, there is an automorphism of algebraic groups $\Psi: \SL_2 (\R) \to \SL_2 (\R)$ and $C \in \SL_2(\R)$ so that $\Gamma _ 2 = \Psi( \Gamma _ 1)$, \ $ u^{(2)}_t = C \Psi( u ^ {(1 )} _ t ) C^{-1}$  so that $\psi$ sends the coset $[g]_{\Gamma _ 1} \in G _ 1 / \Gamma _ 1$ to $ [C\Psi(g)]_{\Gamma_2}$. By measurable isomorphism we mean that $\psi$ is $1$--$1$ onto between conull measurable subsets of $G _ i / \Gamma _ i$, sends $m_1$ to $m_2$ and intertwines with the corresponding $\R$-actions i.e. for every $t \in \R$
\[
\psi (u^{(1)}_t.x) = u^{(2)}_t.\psi(x) \qquad \text{$\mu _ 1$-a.e.}
.\]

This rigidity of measurable isomorphisms was generalized by Witte Morris \cites{witte1985rigidity,witte1987zero}. Rigidity of unipotent flows under measurable isomorphisms can also be deduced from the much stronger measure classification results of Ratner in \cites{ratner1990measure,ratner1991measure} (this is explicitly mentioned by Ratner in \cite[Cor.~6]{ratner1990measure}). We refer the reader to Witte Morris's book \cite{Morris05Ratner} for more details and background.


\medskip

The isomorphism rigidity theorem tells us that a weaker notion of isomorphism between two flows in the class of unipotent flows on quotients of semisimple groups, say $G _ 1 / \Gamma _ 1$ and $G _ 2 / \Gamma _ 2$ implies a much stronger notion of isomorphism. Namely $G_1$ and $G_2$ are isogenous as algebraic groups, and if we chose $G_1$ and $G_2$ (as we may) so that $G_i$ acts faithfully\footnote{By replacing, if needed, the $G_i$ with an appropriate quotient by a finite subgroup in the center of $G_i$.} on $G_i/\Gamma_i$, then in fact $G_1$ and $G_2$ need to be isomorphic with the isomorphism taking $\Gamma_1$ to  $\Gamma_2$ and $u^{(1)}_t$ to a conjugate of $u^{(2)}_t$.

In this paper we allow a much weaker notion of equivalence between flows --- known as \textbf{monotone equivalence} (we shall also use interchangeably the term \textbf{Kakutani equivalence}, though strictly speaking the latter term is most often used for a closely related notion for $\Z$-actions) and investigate if a similar rigidity holds. It turns out that for this equivalence condition a rather striking dichotomy holds: one parameter unipotent flows on quotients of semisimple groups fall into two categories: 
\begin{itemize}
    \item unipotent flows that are loosely Kronecker and hence are all monotone equivalent to each other;
    \item the non loosely Kronecker unipotent flows  where (the weak) monotone equivalence implies the much stronger (algebraic) equivalence as above.
\end{itemize}
We stress that in the loosely Kronecker case monotone equivalence does \textbf{not} imply isomorphism as measure preserving flows. The full details will appear in \cite{Lindenstrauss2022}; we present here some of the ingredients as well as an outline of the argument. An explicit classification of which unipotent flows are loosely Kronecker was obtained by  Kanigowski, Vinhage and the second named author in \cite{kanigowski2021kakutani} (see below).

\medskip

Let $(X _ i, \mathcal B_i,  \mu_i)$ for $i = 1, 2$ be two standard Borel probability measure spaces, and let $u^{(i)}_t : X_i \to X_i$ be an ergodic flow ($\R$-action) on $X _ i$ preserving $m _ i$. The two flows $(X _ i, \mathcal{B}_i, \mu_i, u^{(i)}_t)$ are \textbf{monotone} (or \textbf{Kakutani}) \textbf{equivalent} if there is a 1-1 and onto measurable map $\psi: X' _ 1 \to X' _ 2$ with $\mu_i(X_i \setminus X'_i)=0$ ($X_i' \subseteq X_i$), and so that $\mu_2$ is in the same measure class as $\psi_*\mu_1$,  taking $u^{(1)}_\bullet$-orbits to  $u^{(2)}_\bullet$-orbits \textbf{preserving the order structure}, i.e. so that if $x, u^{(1)}_t.x \in X ' _1$ with $t>0$ then $\psi ( u^{(1)}_t.x) = u^{(2)}_\tau.\psi (x)$ for $\tau > 0$. The condition that $\psi$ preserves the order structure on the orbits is essential, as any two ergodic flows will be orbit equivalent (by work of Ornstein and Weiss, this extends to actions of general amenable groups).

Given an ergodic measure preserving flow $(X , \mathcal B, \mu, u_t)$ and a function $\alpha\in L_+^1(X,\mathcal{B},\mu)$ (one can think of $1/\alpha$ as the time change velocity) we can define a new flow on $X$, monotone equivalent to the original flow, via \textbf{time change}. This new flow $u_{\alpha,t}$ is given by
    \[
    u_{\alpha,\tau}(x)=u_{t}(x) \qquad\text{$\tau$ and $t$ satisfy $\int_0^{t}\alpha(u_sx)ds=\tau$.}
    \]
    The new flow $u_{\alpha,t}$ preserves a measure $\mu^\alpha$ in the same measure class as~$\mu$, with $d \mu^\alpha = \alpha d\mu / (\int\alpha d \mu)$.

\medskip
    
Monotone equivalence can be characterized in terms of time change:  $(X_1 , \mathcal B_1, \mu, u^{(1)}_t)$ is monotone equivalent to $(X_2 , \mathcal B_2, \mu_2, u^{(2)}_t)$ iff there is an $\alpha\in L_+^1(X_1,\mathcal{B}_1,\mu_1)$ so that the time changed system $(X_1 , \mathcal B_1, \mu_1^\alpha, u^{(1)}_{\alpha,t})$ is measurably isomorphic to $(X_2 , \mathcal B_2, \mu_2, u^{(2)}_t)$. 
If $\int\alpha d\mu=1$ we will say that $\psi$ is an \textbf{even Kakutani equivalence}; this implies that \[\psi(u^{(1)}_t x) = u^{(2)}_{t+o_x(t)}\,\psi(x) \qquad\text{a.s.\ as $t \to \infty$}.\]
A monotone equivalence $\psi:(X_1 , \mathcal B_1, \mu, u^{(1)}_t) \to (X_2 , \mathcal B_2, \mu_2, u^{(2)}_t)$ gives rise to a \textbf{cocycle} $\tau(x,t)$ on $(X_1 , \mathcal B_1, \mu, u^{(1)}_t)$ defined by the relation $\psi( u^{(1)}_t.x) = u^{(2)}_{\tau(x,t)}.\psi(x)$. Recall that two cocycles $\tau, \tau'$ on $X_1$ are \textbf{cohomologous} if there is a $w:X_1\to\R$ so that 
\[\tau'(x,t)=\tau(x,t)+w(u^{(1)}_t.x)-w(x).\]

\begin{Definition}
Let $\psi$ and $\psi'$ be two Kakutani equivalences between $(X_1 , \mathcal B_1, \mu, u^{(1)}_t)$ and $(X_2 , \mathcal B_2, \mu_2, u^{(2)}_t)$ inducing cocycles $\tau,\tau'$ on $X_1$ as above. We say $\psi$ is \textbf{cohomologous} to $\psi'$ if there is a $w:X_1\to\R$ 
  \[
    \psi'(x)= u^{(2)}_{w(x)}\circ \psi(x).
    \]
In particular $\tau'(x,t)=\tau(x,t)+w(u^{(1)}_t.x)-w(x)$ hence $\tau$ and $\tau'$ are cohomologous.
\end{Definition}

In the context we study in this paper, i.e.  when both $X_1$ and $X_2$ are quotients of semisimple groups, and $u^{(1)}_t,u^{(2)}_t$ are one parameter unipotent groups, it is easy to ``rectify'' any monotone equivalence between these flows to an even Kakutani equivalence (see Proposition~\ref{prop:evenKak}). Thus we may restrict ourselves to even Kakutani equivalences. The following is our main theorem:

\begin{Theorem}\label{thm:main}
For $i=1,2$, let $G_i$ be the identity component of a real semisimple linear algebraic group without compact factors, $\Gamma_i<G_i$ a lattice, \ $m_i$ the probability measure on $G_i/\Gamma_i$ induced by Haar measure on $G_i$, \ $u_t^{(i)}$ a one-parameter unipotent subgroup of $G_i$, and $\mathfrak{g}_i$ the Lie algebra of $G_i$. Assume for $i=1,2$, the group $u_t^{(i)}$ acts ergodically on $(G_i/\Gamma_i,m_i)$ and that $G_i$ acts faithfully on $G_i/\Gamma_i$. Then $(G_1/\Gamma_1,m_1,u_t^{(1)})$ and $(G_2/\Gamma_2,m_2,u_t^{(2)})$ are Kakutani equivalent if and only if one of the following holds
\begin{enumerate}[label=\textup{(A\arabic*)}]
    \item \label{i:LB} for both $i=1,2$, we have that $\mathfrak{g}_i=\mathfrak{sl}_2(\R)\oplus \mathfrak{g}'_i$ and the generator of $u_t^{(i)}$ is of the form 
    \[                     \begin{pmatrix}
                       0 & 1 \\
                       0 & 0 
                     \end{pmatrix}
                  \times 0 \in \mathfrak{g}_i;\]
   \item There exist an isomorphism $\phi:G_1\to G_2$ and $C\in G_2$ such that $\phi(\Gamma_1)=\Gamma_2$ and $\phi(u^{(1)}_t) = C u^{(2)}_t C^{-1}$. Moreover, any even Kakutani equivalence between these systems is cohomologous to an actual isomorphism.
\end{enumerate}
\end{Theorem}

\medskip

Recall that $(X,\mathcal{B},\mu,T_t)$ is \textbf{loosely Kronecker}\footnote{Such systems are also known as Zero entropy loosely Bernoulli or standard, though the latter term means different things in different contexts.} if $T_t$ is Kakutani equivalent to a linear irrational flow on $\mathbb{T}^2$. Ratner in \cite{ratner1978horocycle} showed that horocycle flows (i.e. the case of $G=\SL_2(\R)$) are loosely Kronecker and thus no rigidity results can hold for $L^1$ time changes of  horocycle flows. Kanigowski, Vinhage
and the second named author \cite[Theorem 1.1]{kanigowski2021kakutani} showed that  $(G/\Gamma,m,u_t)$ in fact is loosely Kronecker if and only if \ref{i:LB} holds. On the other hand Ratner showed that if one assume the time change is done according to a function $\alpha$ that is H\"older in the $\SO_2(\R)$-direction then this implies isomorphism \cite{ratner1986rigidity}.
Further rigidity statements under assumptions e.g.\ like in \cite{ratner1986rigidity} can be found in \cites{tang2020new,afr22}. 

On the other hand, Ratner showed in \cite{ratner1979cartesian} that if $X_1 = \SL_2(\R) / \Gamma_1 \times \SL_2(\R) / \Gamma_1$ and $u^{(1)}_t = \left(\begin{pmatrix}
                       1 & t \\
                       0 & 1 
                     \end{pmatrix},\begin{pmatrix}
                       1 & t \\
                       0 & 1 
                     \end{pmatrix}\right)$
then $(X_1,\mathcal{B},m_1,u^{(1)}_t)$ is not loosely Kronecker hence some rigidity holds --- at the very least, the product of two horocycle flows is not monotone equivalent to a single horocycle flow. Later \cite{ratner1981some} Ratner defined an invariant  under monotone equivalence she called the Kakutani invariant, and used it to show that a product of $k$ horocycle flows is not isomorphic to a product of $\ell$ horocycle flows if $k\neq\ell$. Ratner's invariant was calculated in the generality we consider in this paper by Kanigowski, Vinhage
and the second named author in \cite{kanigowski2021kakutani}. In particular, Ratner's invariant turns out to depend only on the group $G_1$ and the one-parameter unipotent group $u^{(1)}_t$, but not on the lattice $\Gamma_1<G_1$.

In her ICM1994 proceedings paper \cite[p.179]{ratner1995interactions}, Ratner asked the following questions:
\begin{enumerate}
    \item Does there exist a $L^1$-time change between product of horocycle flows on $\operatorname{SL}_2(\R)\times\operatorname{SL}_2(\R)$ with different lattices?
    \item Does the isomorphism rigidity hold for smooth time changes of unipotent flows?
\end{enumerate}

Our main theorem, Theorem \ref{thm:main}, clearly answers the first question, and moreover shows that except for the loosely Kronecker case \textbf{no smoothness assumption} on the time changes is needed. In particular:

\begin{Corollary}\label{cor:rigidity}
Let $G$ be an identity component of a real semisimple algebraic group, $\Gamma< G$ a lattice such that $G$ acts faithfully on $G/\Gamma$, $m$ the probability measure on $G/\Gamma$ induced by Haar measure on $G$ and $u_t$ a one-parameter unipotent subgroup acting ergodically on $(G/\Gamma,m)$. Then $(G/\Gamma, m,u_t)$ has time change rigidity if and only it is not loosely Kronecker.
\end{Corollary}

In a recent paper, Gerber and Kunde \cite{gerber2021anti} investigated the descriptive set theoretic complexity of the general monotone equivalence problem for flows, showing that in full generality one cannot classify measure preserving flows up to monotone equivalence using any single (or indeed, countably many) invariants. In their paper they ask about the special case of unipotent flows and Ratner's Kakutani invariant. Our result provides an answer for this question:
\begin{Corollary}
Let $G=\operatorname{SL}_2(\R)\times\operatorname{SL}_2(\R)$, $\Gamma_1,\Gamma_2< G$ two non-conjugate lattices, $u_t^{(i)}=\left(\begin{pmatrix}
                       1 & t \\
                       0 & 1 
                     \end{pmatrix},\begin{pmatrix}
                       1 & t \\
                       0 & 1 
                     \end{pmatrix}\right)$ and $m_i$ the probability measure on $G/\Gamma_i$ induced by Haar measure on $G$ for $i=1,2$. Then $(G/\Gamma_1,m_1,u_t^{(1)})$ and $(G/\Gamma_2,m_2,u_t^{(2)})$ are not Kakutani equivalent while they have the same value of Kakutani invariant.
\end{Corollary}

\subsection*{Acknowledgement}
The authors would like to thank Benjamin Weiss for encouragement and helpful discussions. D.W. would also like to thank Svetlana Katok for her encouragements and Adam Kanigowski, Kurt Vinhage and Andreas Wieser for many helpful discussions and comments. We are grateful to the anonymous referee's helpful comments and suggestions.

\section{Some preliminaries}
In this section we present some notations that are needed to state the key steps in the proof of Theorem \ref{thm:main}.
\subsection{\texorpdfstring{$(\delta,\epsilon,R)$}--two sides matchable}\label{sec:twoSideM}
In analogy to Feldman's $\bar f$ metric we will make use of the notion of $(\epsilon,\delta,R)$-two sided matching between two points. For more details and background we refer the reader to \cites{feldman1976new,katok1977monotone,ornstein1982equivalence,ratner1981some}.

Let $l$ denotes the Lebesgue measure on $\R$ and consider for $R>1$ the orbit segment
$I_R(x)=\{T_sx:s\in[-R,R]\}$.

\begin{Definition}[cf.~{\cite[Def.~1]{ratner1981some}}]
Let $T_t$ be a flow on $(X,\mathcal{B},\mu)$ and $d$ a metric on $X$. Suppose $\delta,\epsilon\in(0,1)$ and $R>1$. We say $x,y\in X$ are \textbf{$(\delta,\epsilon,R)$-two sides matchable} if there exist a subset $A=A_{(x,y)}\subset[-R,R]$ with $l(A)>(1-\epsilon)2R$, and an increasing absolutely continuous map $h=h_{(x,y)}:A \to [-R,R]$, such that
\begin{enumerate}
    \item $d(T_{h(t)}x,T_ty)<\delta$ for all $t\in A$;
    \item $h(0)=0$ and $0\in A$;
    \item $|h'(t)-1|<\epsilon$ for all $t\in A$.
\end{enumerate}
We call $h$ an \textbf{$(\delta,\epsilon,R)$-two sided matching} between $I_R(x)$ to $I_R(y)$.
\end{Definition}


\subsection{Our setup}
In this paper, we take $G$ to be the identity component of a real semisimple linear algebraic group. Cf.\ e.g. \cite[Ch.~1]{MargulisDiscrete1991} for a concise introduction to the theory of such groups.
Let $\mathfrak{g}$ denote the Lie algebra of $G$.

Let $d_G$ denote a right invariant Riemannian metric on $G$, which also induces a Riemannian metric $d_{G/\Gamma}$ on $G/\Gamma$. Moreover, the corresponding Riemannian volume on $G$ defines a Haar measure $\tilde{m}$ on $G$ and thus induces a probability measure $m$ on $G/\Gamma$.

\subsection{\texorpdfstring{$\mathfrak{sl}_2(\R)$}{sl2(R)}---triple, chain basis and optimal matching function}\label{sec:sl2triple}
In order to give a precise description of the divergence rate for nearby points in $G/\Gamma$, the following special basis in Lie algebra $\mathfrak{g}$ is very helpful. 
We make us of the standard notations about exponential map, adjoint action and their relations; we refer the reader to \cite{Knapp96Lie} for more details.

Let $\bu\in\mathfrak{g}$ be a nilpotent element; by Jacobson-Morozov theorem \cite{jacobson1979lie}, there exists a homomorphism $\varphi:\mathfrak{sl}_2(\R)\to\mathfrak{g}$ such that
\begin{equation*}
    \begin{aligned}
        \varphi\left(\left(
                     \begin{array}{cc}
                       0 & 1 \\
                       0 & 0 \\
                     \end{array}
                   \right)\right)=\bu, \ \ \varphi\left(\left(
                     \begin{array}{cc}
                       1 & 0 \\
                       0 & -1 \\
                     \end{array}
                   \right)\right)=\ba, \ \ \varphi\left(\left(
                     \begin{array}{cc}
                       0 & 0 \\
                       1 & 0 \\
                     \end{array}
                   \right)\right)=\bu^-,
    \end{aligned}
\end{equation*}
where $\ba\in\mathfrak{g}$ is an $\R$-diagonalizable element and $\bu^-\in\mathfrak{g}$ is a nilpotent element. We shall use the following three one-parameter subgroups:
\begin{itemize}
    \item the one-parameter unipotent subgroup 
    \[\displaystyle U=\{u_t=\exp(t\bu):t\in\R\};\]
    \item the corresponding one-parameter diagonalizable subgroup 
    \[\displaystyle A=\{a_t=\exp(t\ba):t\in\R\};\]
    \item the ``opposite'' one-parameter subgroup 
    \[\displaystyle V=\{v_t=\exp(t\bu^-):t\in\R\}.\]
\end{itemize}
Let $\mathfrak{h}=\operatorname{span}\{\bu,\ba,\bu^-\}$ and let $H< G$ be the connected Lie subgroup (in the Hausdorff topology) generated by $\exp(\mathfrak{h})$.

\medskip

By the well-known classification of finite dimensional Lie representations of $\mathfrak{sl}_2(\R)$, we can find a basis of $\mathfrak{g}$ of the form
\begin{equation}\label{eq:chainbasis}
    \{\bu,\ba,\bu^-,\bx^{0,1},\ldots,\bx^{m_1,1},\ldots,\bx^{0,n},\ldots,\bx^{m_n,n}\},
\end{equation}
where for each $1\leq j \leq n$, we have that $\bx^{0,j},\ldots,\bx^{m_j,j}$ is an irreducible representation of $\mathfrak{h}$ under $\ad$ with $\ad_{\bu}\bx^{i,j}=\bx_{i+1,j}$.
We call this basis the chain basis for $\mathfrak{g}$ with respect to $\mathfrak{h}$.

Given a subset $W\subset\mathfrak{g}$, we define $C_{\mathfrak{g}}(W)$ to be the common centralizer of all $\bw\in W$, i.e.
\begin{equation}
    C_{\mathfrak{g}}(W)=\{\by\in\mathfrak{g}:\ad_{\bw}(\by)=0\ \ \forall\bw\in W\}.
\end{equation}

\medskip

In order to find the best matching between the $U$-orbits of two nearby points, say $u_t.x$, and $u_t.wx$ one needs to modify the element of $U$ used in one of these two points to compensate for the shearing behaviour of the unipotent flow. Explicitly, we have the following:

\begin{Lemma}\label{lem:matchingFunction}
Suppose
\[
w=\exp(\vartheta_{\bu^-}\bu^-)\exp(\vartheta_{\ba}\ba)\in G.
\]
Then there exists a rational function $\phi(t)$ (that depends on $w$) such that 
\begin{equation}\label{eq: 2by2matrices calulation}
\exp(\phi(t)\bu)w\exp(-t\bu)=\exp(\vartheta_{\bu^-,t}\bu^-)\exp(\vartheta_{\ba,t}\ba),
\end{equation}
so that if $\absolute{t} < \epsilon \absolute{\vartheta_{\bu^-}^{-1}}$ and  $\absolute{\vartheta_{\ba}}<\epsilon$ then 
\[
\absolute{\vartheta_{\ba,t}}<2\epsilon,\ \ \ \ \absolute{\vartheta_{\bu^-,t}}<2\absolute{\vartheta_{\bu^-}}.
\]
Moreover, $\absolute{\phi'(t)-1}<\sqrt{\epsilon}$.

\end{Lemma}

\medskip

This lemma can be viewed as a special case of what Ratner calls the H-property (see \cite[Def. 1]{ratner1983horocycle}). Equation \eqref{eq: 2by2matrices calulation} and the following inequalities can be obtained by direct computation of $2\times2$ matrices. Cf. \cite[Lemma 5.3]{kanigowski2021kakutani} for more details.

\subsection{Kakutani-Bowen ball}\label{sec:Kakutanibowen}


We can represent elements $g\in G$ sufficiently close to identity using the following local chart based on our choice of basis for $\mathfrak{g}$ in \eqref{eq:chainbasis}:
\begin{equation}\label{eq:gDecom}
    g=\exp(\vartheta_{\bu}(g)\bu)\exp(\vartheta_{\bu^-}(g)\bu^-)\exp(\vartheta_{\ba}(g)\ba)\breve g,
\end{equation}
where
\begin{equation*}
    \begin{aligned}
    \breve g&=\exp(\sum_{j=1}^n\sum_{i=0}^{m_j}\vartheta_{i,j}(g)\bx^{i,j}).
    \end{aligned}
\end{equation*}

\begin{Definition}
\begin{enumerate}
\item
For $\epsilon>0$ and $R>1$, let $\operatorname{Bow}(R,\epsilon)$ be the set
\begin{equation*}
  \{g\in G,d_G(\exp(t\bu)g\exp(-t\bu),e)<\epsilon\textup{ for }\forall t\in[-R,R]\},
\end{equation*}
and define for $x \in G/\Gamma$ the \textbf{$(R,\epsilon)$-Bowen ball} around $x$ to be
\[
\operatorname{Bow}(R,\epsilon,x )=\operatorname{Bow}(R,\epsilon).x.
\]
\item 
For $\epsilon>0$ and $R>1$, let $\operatorname{Kak}(R,\epsilon)$ be the set of $g \in G$ satisfying the following:
\begin{enumerate}
    \item\label{item:Kak1} $\absolute{\vartheta_{\bu^-}(g)}<\frac{\epsilon}{R}$,
    \item\label{item:Kak2} $\absolute{\vartheta_{\ba}(g)}<\epsilon$,
    \item\label{item:Kak3} $\absolute{\vartheta_{\bu}(g)}<\epsilon$,
    \item\label{item:Kak4} $\breve g\in\operatorname{Bow}(R,\epsilon)$,
\end{enumerate}
where $\vartheta_{\bu^-}$, $\vartheta_{\ba}$, $\vartheta_{\bu}$ and $\breve g$ as in \eqref{eq:gDecom}. For $x\in G/\Gamma$, we define the \textbf{$(R,\epsilon)$-Kakutani-Bowen ball} around it to be \[\Kak(R,\epsilon,x)=\Kak(R,\epsilon).x,\] and similarly set $\Kak(R,\epsilon,\tilde{x})=\Kak(R,\epsilon).\tilde x$ for $\tilde x \in G$.
\end{enumerate}
\end{Definition}

It is easy to see using the ``time change'' given in Lemma~\ref{lem:matchingFunction} that if $x \in G/\Gamma$ and $y \in \Kak(R,\epsilon,x)$ then $x$ and $y$ are $(C\epsilon,\sqrt{\epsilon},R)$-matchable for an appropriate constant $C$, and similarly for $\tilde x,\tilde y \in G$.
On $G$, the converse also essentially holds, and moreover the best matching among all two sided matching is essentially given by the matching defined in \S\ref{sec:twoSideM}. This turns out to be also true in $G/\Gamma$, but this is a much more delicate statement and only holds \textbf{if} the action of $u_t$ on $G/\Gamma$ is not loosely Kronecker --- indeed, this fact is essentially equivalent to our Main Lemma presented in \S\ref{sec:mainlemma}.

\begin{Proposition}
There exist constants $\Cnew\label{matchKak}>0$ and $L$
such that for any $\epsilon,\delta>0$ small enough the following holds.
Suppose $\tilde{x}\in G$, $d_G(g,e)<\delta$ and $\tilde{y}=g.\tilde{x}$. Let $h:\R\to\R$ be a $C^{\infty}$ function such that $h(0)=0$ and $\absolute{h'(t)-1}<\epsilon$. Then we can cover the set  
\[
I_{\tilde{x},\tilde{y}}=\{t\geq 0:d_{G}(\exp(h(t)\bu)\tilde{x},\exp(t\bu)\tilde{y})<4\delta\} 
,\] 
by $k \leq L$ disjoint closed intervals $[b_i,d_i]$, so that on each interval
\begin{equation*}
 d_G(\exp((\phi(t)-\phi(b_i)+h(b_i))\bu)\tilde{x},\exp(t\bu)\tilde{y})<\Cold{matchKak}\delta\qquad \forall t\in [b_i,d_i]
\end{equation*}
with
$\phi(t)$ is as in Lemma \ref{lem:matchingFunction}. Moreover for all $i=1,\dots,k$ and $t \in [b_i,d_i]$
\[
\exp((\phi(t)-\phi(b_i)+h(b_i))\bu)\tilde{x}\in\Kak(d_i-b_i,\Cold{matchKak}\delta,\exp(t\bu)\tilde{y}).
\]
\end{Proposition}

\subsection{Two relations}\label{sec:twoRelations}
Let $K\subset G/\Gamma$ be a compact set, $x,y\in K$ and $\epsilon$ small enough so that for any $x \in K$, the map $g \mapsto g.x$ is injective on the $\epsilon$-ball around $e \in G$. Suppose that $y=g.x$ for $g \in G$ with $d_{G}(g,e)<\epsilon$ and moreover that $d_{G/\Gamma}(u_t.x,u_s.y)<\epsilon$ for some $t,s\geq0$. Lift $x$ to a point $\tilde x$ in $G$ and set $\tilde y=g\tilde x$, so that $d_G(\tilde x,\tilde y)<\epsilon$. Then there exists a unique $\gamma\in\Gamma$ such that $d_G(\exp(t\bu)\tilde x,\exp(s\bu)\tilde y\gamma)<\epsilon$. Note that the conjugacy class $[\gamma]$ of $\gamma$ does not depend on the lift $\tilde x$ we choose for $x$.
If $x,y,s,t,[\gamma]$ satisfy the above we shall write
\[(x,y)\overset{[\gamma]}{\leadsto}(u_t.x,u_s.y);
\]
we also write $(x,y)\overset{e}{\leadsto}(u_t.x,u_s.y)$ for $(x,y)\overset{[e]}{\leadsto}(u_t.x,u_s.y)$ and 
\[(x,y)\overset{\neq e}{\leadsto}(u_t.x,u_s.y)\]
if $(x,y)\overset{[\gamma]}{\leadsto}(u_t.x,u_s.y)$ for $\gamma\neq e$.

\section{Reduction to a good Kakutani equivalence}
We begin with two reductions for any general Kakutani equivalence between the actions of unipotent one-parameter subgroups. The second step of these reductions also holds for any general abstract ergodic systems.

For $i=1,2$, we assume that $G_i$ is the identity component of a real semisimple linear algebraic group without compact factors, $\Gamma_i< G_i$ a lattice and $m_i$ a probability measure on $G_i/\Gamma_i$ induced from Haar measure on $G_i$. Let $u_t^{(i)}$ be a one-parameter unipotent subgroup acting ergodically on $(G_i/\Gamma_i,m_i)$ for $i=1,2$ and $\psi$ an arbitrary Kakutani equivalence between $(G_1/\Gamma_1,u_t^{(1)},m_1)$ and $(G_2/\Gamma_2,u_t^{(2)},m_2)$.  

The first step of the reduction is using the diagonal subgroup to normalize our one-parameter unipotent subgroup, which reduces $\psi$ to an even Kakutani equivalence:
\begin{Proposition}\label{prop:evenKak}
There exists $s_0\in\R$ such that $a_{s_0}^{(2)}\circ\psi$ is an even Kakutani equivalence, where $a_s^{(2)}$ is the diagonalizable subgroup arising from $\mathfrak{sl}_2(\R)$-triple for the generator of $u_t^{(2)}$ in \S\ref{sec:sl2triple}.
\end{Proposition}

The target of the second step of our reduction is to obtain a good control of the time change function $\alpha$. In order to simplify the notation, we introduce the following definition:
\begin{Definition}\label{def:goodevenkak}
Let $T_t$ and $S_t$ be two ergodic flows acting on $(X,\mathcal{B},\mu)$ and $(Y,\mathcal{C},\nu)$ respectively. For any $\epsilon>0$, we say $\psi$ is an \textbf{$\epsilon$-controllable Kakutani equivalence} between $T_t$ and $S_t$ if the following holds:
\begin{enumerate}
    \item $\psi$ is an even Kakutani equivalence between $T_t$ and $S_t$ with time change function $\alpha$;
    \item $\operatorname{esssup}\absolute{\alpha-1}<\epsilon$;
    \item there exists a full measure set $K\subset X$ such that for every $x\in K$, $\alpha(T_tx)$ is a $C^{\infty}$ function in $t$.
\end{enumerate}
\end{Definition}



\begin{Lemma}\label{lem:goodevenkak}
Let $\psi$ be an even Kakutani equivalence between two ergodic flows. For any given $\epsilon>0$, there exists an $\epsilon$-controllable Kakutani equivalence $\tilde{\psi}$ that is cohomologous to $\psi$.
\end{Lemma}

This lemma (which is a lemma about general ergodic systems and does not use any special properties of unipotent flows) follows from a combination of {\cite[Proposition 2.3]{katok1977monotone}} and {\cite[Theorem 1.4]{ornstein1982equivalence}}. 

As a result of these two reductions, it is permissible to assume that $\psi$ is a $\enew\label{eKak}$-controllable Kakutani equivalence, with $\eold{eKak}>0$ a small fixed constant, an assumption we make from this point on.

\section{Main lemma}\label{sec:mainlemma}
One of the key ingredients of the proof of Theorem~\ref{thm:main} is the following proposition, which essentially shows that Kakutani-Bowen balls are preserved under Kakutani equivalence map. This type of proposition first appeared in \cite{ratner1979cartesian} for Cartesian products of horocycle flows on compact quotients of $\SL_2(\R)$. In order to calculate an invariant for Kakutani equivalences introduced by Ratner \cite{ratner1981some} for more general one-parameter unipotent flows, Kanigowski, Vinhage and the second named author extended this analysis to general semisimple linear Lie groups in \cite{kanigowski2021kakutani}, though for our purposes we need a sharper form of this important result.

\begin{Lemma}[Main Lemma]\label{lem:main}
There exist $\denew\label{027delm1},\Rnew\label{027Rlm1}>0$ and a compact set $\Knew\label{028Klm1}\subset G_1/\Gamma_1$ with $m_1(\Kold{028Klm1})>0.99$ such that for any $\delta\in(0,\deold{027delm1})$, there exists $\epsilon\in(0,\delta)$ such that if $x,y\in\Kold{028Klm1}$ and $R>\Rold{027Rlm1}$ satisfy
\[
x\in\Kak(R,\epsilon,u_{t_0}^{(1)}.y)
\]
for some $\absolute{t_0}<\epsilon R$, then we have
\[
\psi(u_{t_0}^{(1)}.y)\in \Kak(R,\delta,u_{t_R}^{(2)}.\psi(x)),
\]
for some $t_R$ satisfying $\absolute{t_R}\leq40\eold{eKak} R$.
\end{Lemma}

\medskip


This main lemma is at the heart of our whole argument as it provides an essential tool to establish that images of two nearby points in a good set under a Kakutani equivalence remain close; a key point is that for this purpose closeness needs to be measured in terms of Kakutani-Bowen balls.

\medskip

In order to establish the main lemma, we first note that as in \S\ref{sec:Kakutanibowen}, $x\in\Kak(R,\epsilon,u_{t_0}^{(1)}.y)$ implies that $x$ and $u_{t_0}^{(1)}.y$ are $(C\epsilon,\sqrt{\epsilon},R)$-two sides matchable for some appropriate constant $C$. On a set of large measure, $\psi$ is continuous, so combining the pointwise ergodic theorem with our assumption that $\psi$ is an $\eold{eKak}$-controllable Kakutani equivalence (cf.\ the end of the previous section)  give that if $x$ and $y$ are in a ``good'' set of large measure and $R$ large enough then $\psi(x)$ and $\psi(u_{t_0}^{(1)}.y)$ are $(\delta,2\eold{eKak},R)$-two sides matchable for some appropriate constant $\delta$.
\medskip

In view of this, Lemma \ref{lem:main} can be reduced to the following lemma, which establishes the connection between two sided matching in $G_2/\Gamma_2$ and Kakutani-Bowen balls. 
\begin{Lemma}\label{lem:longHamming}
There is a compact subset $\Knew\label{029Klh1}\subset G_2/\Gamma_2$ with $m_2(\Kold{029Klh1})\geq0.99$ and constants $\Cnew\label{027Clh1},L$ such that for any $\epsilon,\delta$ sufficiently small and $R$ sufficiently big the following holds. Let $x \in G_2/\Gamma_2$ and $ y \in \Kold{029Klh1}$ be $(\Cold{027Clh1}\delta,\epsilon,R)$-two sides matchable with matching function $h$.
Then there exist lifts $\tilde{x}, \tilde{y} \in G_2$ of $x,y$, \ $\gamma_R\in\Gamma_2$, \ $\absolute{s_R}\leq R$ and $R'\in\left[\frac{R}{L}, R\right]$ such that for both  $t=0$ and $R'$,
\begin{equation*}
    \exp(\left(h(s_R+t)\right)\bu_2)\tilde{x}\in\Kak(R',\delta)\exp((s_R+t)\bu_2)\tilde{y}\gamma_R.
\end{equation*}
\end{Lemma}

\medskip
Roughly speaking, the above lemma says that if two points are two sides matchable for sufficiently long time, then they also stay close in the universal cover for a long time of the same order of magnitude. 


\medskip
Lemma \ref{lem:longHamming} can be derived from the following lemma, where relation $\overset{\neq e}{\leadsto}$ is defined in \S\ref{sec:twoRelations}. In order to simplify the notation, we define $x_t=u_t^{(2)}.x$ and $y_t=u_t^{(2)}.y$ for $x,y\in G_2/\Gamma_2$.
\begin{Lemma}\label{Lem:expgap}
There exist $\Cnew\label{Cexp},w,\Rnew\label{Rexp}>0$ and a compact set $\Knew\label{Kexp}\subset G_2/\Gamma_2$ with $m_2(\Kold{Kexp})>0.99$ such that for every $R\geq\Rold{Rexp}$ the following holds. Suppose that $y\in\Kold{Kexp}$, $x\in\Kak(R,\delta,y)$, $x_{h(s)}\in\Kak(R,\delta,y_s)$ and $(x,y)\overset{\neq e }{\leadsto}(x_{h(s)},y_s)$, then we have, 
\[\min(s,h(s))\geq\Cold{Cexp}R^{1+w}.\]
\end{Lemma}

Such a result first appeared in Ratner \cite{ratner1979cartesian} for the Cartesian product of horocycle flows. By using techniques in algebraic group theory, we generalize it to arbitrary one-parameter unipotent subgroups. Together with some combinatorial estimates, Lemma \ref{Lem:expgap} shows that for most points, the two sided matching on $G_2/\Gamma_2$ implies the two sided matching in the same time scale on $G_2$ for their corresponding lifts. This will give the proof of Lemma \ref{lem:longHamming} and hence also Lemma \ref{lem:main}.



\section{Compatibility of diagonalizable group}
For any $\epsilon, \delta$, if $\eta$ is sufficiently small, $x\in G_1/\Gamma_1$ then $x$ and $a^{(1)}_{\eta} x$ are $(\delta,\epsilon,R)$ two sides matchable for any $R$.
Using Lemma~\ref{lem:main} (the ``Main Lemma'') this implies that $\psi(a^{(1)}_{\eta}x) \in N_{G_2}(U^{(2)}) \psi(x)$ with $N_{G_2}(U^{(2)})$ denoting the normalizer in $G_2$ of $U^{(2)}$, where $U^{(2)}$ is the one-parameter unipotent subgroup defined in \S\ref{sec:sl2triple}. This can further be upgraded to the following statement:
\begin{Proposition}\label{prop:geoCom}
There is an element $c$ in the centralizer of $U^{(2)}$ so that if $\psi'(x)=c\psi(x)$ then a.e.
\[
\psi'(a^{(1)}_{1}.x)\in U^{(2)} a_{1}^{(2)}.\psi'(x).
\]
\end{Proposition}

Without loss of generality we may assume $\psi'=\psi$, and then Proposition \ref{prop:geoCom} says that there is a measurable function $t(x):X_1\to\R$ so that for a.e.~$x$
\begin{equation}\label{eq:geo1com}
  \psi(a_1^{(1)}.x)=u_{t(x)}^{(2)}a_1^{(2)}.\psi(x).
\end{equation}
Iterating, we can get
\begin{equation}\label{eq:geoIter}
 \psi(a_n^{(1)}.x)=u_{e^{2n-2}t(x)+e^{2n-4}t(a_1^{(1)}.x) + \dots +t(a_{n-1}^{(1)}.x)}^{(2)}a_n^{(2)}.\psi(x).
\end{equation}

Unfortunately, we have no information on $t(x)$ (we certainly do \textbf{not} know it is in $L^1$) so $e^{2n-2}t(x)+e^{2n-4}t(a_1^{(1)}.x) + \dots +t(a_{n-1}^{(1)}.x)$ could be very large.

We overcome this significant hurdle by perturbing the points $a_j^{(1)}.x$ into a set where $t(\cdot)$ is bounded, using a certain
 $L^2$ ergodic theorem for $\SL_2(\mathbb R)$-actions described in the next section.

In her paper \cite{ratner1986rigidity}, Ratner used a statement analogous to Proposition~\ref{prop:geoCom} to show that \textbf{if the time change is sufficiently smooth}, a Kakutani equivalence between two horocycle flows has to come from an isomorphism. We follow a similar strategy here, using 
equation \eqref{eq:geo1com} as a crucial ingredient. To do that,  we consider on the ``conjugation'' of $\psi$ under the diagonalizable groups $a_t^{(i)}$
\begin{equation}\label{eq:psiConA}
  \psi_{n}(x)=a_{-n}^{(2)}\psi\left(a_{n}^{(1)}.x\right), \ \ x\in G_1/\Gamma_1.
\end{equation}
Combining the even Kakutani equivalence assumption and pointwise ergodic theorem, the limit of $\psi_n$ (\textbf{if} it exists) will be a measurable isomorphism between the corresponding unipotent flows and thus the only remaining task is the existence of the limit of $\psi_n$. In \cite{ratner1986rigidity}, the regularity assumption of the time changes is used to establish the existence of the limit of $\psi_n$ along a subsequence, thus establishing time change rigidity in this case were the time change is assumed to have some apriori regularity. A similar strategy has also been used very recently by Artigiani, Flaminio and Ravotti \cite{afr22}, however they need to \textbf{assume} convergence of the $\psi_n$ to get an isomorphism.
Due to the lack of regularity assumption in our situations, significant changes are needed. 

\section{Spectral gap and an \texorpdfstring{$\SL_2(\mathbb{R})$}{SL2(R)}-ergodic theorem}
The aim of this section is to state a $\SL_2(\R)$-pointwise ergodic theorem, which enable us to find good points in pushforward of small $\epsilon$-balls in $H_1$ under $a_n^{(1)}$. We remark that the spectral gap for $L^2(G_1/\Gamma_1)$ is only used in this section to obtain the necessary pointwise behaviour.

Our $\SL_2(\R)$-ergodic theorem is following:
\begin{Theorem}\label{thm:sl2ergodic}
Let $\rho$ be a unitary representation of $H=\SL_2(\R)$ on a Hilbert space $\mathcal H_{\rho}$ with no fixed vectors and a spectral gap. Then there exists $\enew\label{eSpectral}>0$ such that for any $\epsilon\in(0,\eold{eSpectral})$ there is a $C(\epsilon)$ so that  
\[
\sum_{n=1}^{+\infty}\norm{\frac{1}{m_{H}(B_{\epsilon}^{H,\norm{\cdot}})}\int_{B_{\epsilon}^{H,\norm{\cdot}}}\rho( ha_{-n})v \,dm_H(h)}^{2}<C(\epsilon)\|v\|^2
\]
for any $v \in \mathcal H_{\rho}$,
where $m_H$ is the Haar measure on $H$ and $B_{\epsilon}^{H,\norm{\cdot}}$ is the $\epsilon$-ball in $H$ around identity with respect to the Hilbert-Schmidt norm.
\end{Theorem}

Recall that a unitary representation $(\rho, \mathcal H)$ of $H$ is said to have a spectral gap if there is some compactly supported probability measure $\nu$ on $H$ so that on the orthogonal complement of the $H$-fixed vectors the norm of the operator $\int \rho(g) d\nu(g)$ is $<1$.


A direct corollary of Theorem \ref{thm:sl2ergodic} is:
\begin{Corollary}\label{cor:sl2erg}
Let $G$ be the identity component of a real semisimple algebraic group without compact factors, $\Gamma< G$ a lattice with $m$ the probability measure on $G/\Gamma$ induced by Haar measure on $G$. Let $H < G$ be locally isomorphic to $\SL_2$  and let $a_t \in H$ be a one-parameter diagonalizable subgroup as in \S\ref{sec:sl2triple}.
Given $\epsilon\in(0,\eold{eSpectral})$ and $\eta\in[0,1)$, then for any $f\in L^2(G/\Gamma,m)$, $m$-a.e. $x\in G/\Gamma$, we have
\[
\lim_{n\to+\infty}\frac{1}{m_H(B_{\epsilon}^{H,\norm{\cdot}})}\int_{B_{\epsilon}^{H,\norm{\cdot}}}f(a_{n(1-\eta)}ha_{n\eta}.x)dm_H(h)=\int_{G/\Gamma}fdm.
\]
\end{Corollary}

\medskip

The deduction of Corollary~\ref{cor:sl2erg} from Theorem~\ref{thm:sl2ergodic} uses the fact that if~$G$ is a semisimple group without compact factors, $\Gamma$ an irreducible lattice, and $H<G$ is as above, then the representation of $H$ on $L^2(G/\Gamma)$ arising from left translations has a spectral gap. More generally, if $G=\prod_i G_i$ and $\Gamma=\prod_i \Gamma_i$ with $\Gamma_i <G_i$ irreducible lattices, and if $H$ projects non-trivially into each $G_i$, then the action of $H$ on $L^2(G/\Gamma)$ has a spectral gap. C.f.\ e.g.~\cites{kleinbock1999log,KelSar09} for more information and background.

\section{Existence of limiting map}

Finally we show that $\psi_n(x)$ converges in an appropriate sense to a limiting map $\varphi(x)$. This limiting map will turn out to be not just an (even) Kakutani equivalence but in fact a measurable isomorphism between the two unipotent flows given by the $u^{(i)} _ t$ action on $G_i/\Gamma_i$.

\begin{Theorem}\label{thm:reno}
For $m_1$-a.e. $x\in G_1/\Gamma_1$, there is a subsequence $\{n_i\}_{i\in\N}$ of natural numbers with full density such that
$\varphi(x)=\lim_{i\to\infty} \psi_{n_i}(x)$ exists and
\begin{enumerate}[label=(B\arabic*)]
    \item\label{item:phiM} $\varphi:G_1/\Gamma_1\to G_2/\Gamma_2$ is a measurable map;
    \item\label{item:Uorbit} $\varphi(x)\in U^{(2)}(\psi(x))$;
    \item\label{item:uPre} $\varphi(u_t^{(1)}.x)=u_t^{(2)}.\varphi(x)$ for $m_1$-a.e. $x\in G_1/\Gamma_1$ and any $t\in\mathbb{R}$;
    \item\label{item:geoPre} $\varphi(a_k^{(1)}.x)=a_k^{(2)}.\varphi(x)$ for $m_1$-a.e. $x\in G_1/\Gamma_1$ and any $k\in\mathbb{N}$.
\end{enumerate}
\end{Theorem}

The challenging part is establishing that there is a subsequence $\{n_i\}_{i\in\N}$ of natural numbers with full density such that
\[
\varphi(x)=\lim_{i\to\infty} \psi_{n_i}(x)
\] exists; the rest follows by standard arguments.

The key observation to prove Theorem \ref{thm:reno} is to consider the iterations $\{a_{n}^{(1)}h_n.x\}_{n\in\N}$ for some $h_n\in B_{\epsilon}^{H_1,\norm{\cdot}}$ instead of $\{a_n^{(1)}.x\}_{n\in\N}$. By using Corollary \ref{cor:sl2erg}, we can always guarantee that $a_{n}^{(1)}h_n.x$ stays in some good sets, i.e. $\absolute{t(a_{n}^{(1)}h_n.x)}$ is bounded by a positive constant, for every $n$. Note that this was not possible for $a_n^{(1)}x$ due to the existence of bad iterations as guaranteed by pointwise ergodic theorem. In order to control the error terms in opposite horocycle directions, i.e. $\exp(\R\bu_2^-)$, we also need to consider iterations $\{a_{n(1-\eta)}^{(1)}h_na_{n\eta}^{(1)}.x\}_{n\in\N}$. With the help of estimates in  representation theory and Lemma \ref{lem:main}, we are able to combine two iterations and obtain controls for all directions. 

\medskip 

Once we obtain Theorem \ref{thm:reno}, the proof of Theorem \ref{thm:main} follows from the combination of isomorphism rigidity theorem for unipotent flows and Ratner's argument in \cite{ratner1986rigidity}. More precisely, \ref{item:phiM}, \ref{item:uPre} together with isomorphism rigidity theorem for unipotent flows give that $G_1$ and $G_2$ are isomorphic, $\Gamma_1$ and $\Gamma_2$ are conjugate up to group isomorphism. Combining \ref{item:Uorbit}, \ref{item:geoPre} and Ratner's argument in \cite[Proof of Theorem $1$]{ratner1986rigidity}, we obtain that $\psi$ indeed is cohomologous to a group isomorphism. This completes the proof.

\bibliographystyle{plain} 
\bibliography{refs} 

\begin{bibdiv}
\begin{biblist}

\bib{afr22}{article}{
      author={Artigiani, M.},
      author={Flaminio, L.},
      author={Ravotti, D.},
       title={On rigidity properties of time-changes of unipotent flows},
        date={2022},
     journal={arXiv preprint arxiv 2209.01253},
}

\bib{feldman1976new}{article}{
      author={Feldman, J.},
       title={New {$K$}-automorphisms and a problem of {K}akutani},
        date={1976},
        ISSN={0021-2172},
     journal={Israel J. Math.},
      volume={24},
      number={1},
       pages={16\ndash 38},
         url={https://doi.org/10.1007/BF02761426},
      review={\MR{409763}},
}

\bib{gerber2021anti}{article}{
      author={Gerber, Marlies},
      author={Kunde, Philipp},
       title={Anti-classification results for the kakutani equivalence
  relation},
        date={2021},
     journal={arXiv preprint arXiv:2109.06086},
}

\bib{jacobson1979lie}{book}{
      author={Jacobson, Nathan},
       title={Lie algebras},
   publisher={Dover Publications, Inc., New York},
        date={1979},
        ISBN={0-486-63832-4},
        note={Republication of the 1962 original},
      review={\MR{559927}},
}

\bib{katok1977monotone}{article}{
      author={Katok, A.~B.},
       title={Monotone equivalence in ergodic theory},
        date={1977},
        ISSN={0373-2436},
     journal={Izv. Akad. Nauk SSSR Ser. Mat.},
      volume={41},
      number={1},
       pages={104\ndash 157, 231},
      review={\MR{0442195}},
}

\bib{kleinbock1999log}{article}{
      author={Kleinbock, D.~Y.},
      author={Margulis, G.~A.},
       title={Logarithm laws for flows on homogeneous spaces},
        date={1999},
        ISSN={0020-9910},
     journal={Invent. Math.},
      volume={138},
      number={3},
       pages={451\ndash 494},
         url={https://doi.org/10.1007/s002220050350},
      review={\MR{1719827}},
}

\bib{Knapp96Lie}{book}{
      author={Knapp, Anthony~W.},
       title={Lie groups beyond an introduction},
      series={Progress in Mathematics},
   publisher={Birkh\"{a}user Boston, Inc., Boston, MA},
        date={1996},
      volume={140},
        ISBN={0-8176-3926-8},
         url={https://doi.org/10.1007/978-1-4757-2453-0},
      review={\MR{1399083}},
}

\bib{KelSar09}{article}{
      author={Kelmer, Dubi},
      author={Sarnak, Peter},
       title={Strong spectral gaps for compact quotients of products of {${\rm
  PSL}(2,\Bbb R)$}},
        date={2009},
        ISSN={1435-9855},
     journal={J. Eur. Math. Soc. (JEMS)},
      volume={11},
      number={2},
       pages={283\ndash 313},
         url={https://doi.org/10.4171/JEMS/151},
      review={\MR{2486935}},
}

\bib{kanigowski2021kakutani}{article}{
      author={Kanigowski, Adam},
      author={Vinhage, Kurt},
      author={Wei, Daren},
       title={Kakutani equivalence of unipotent flows},
        date={2021},
        ISSN={0012-7094},
     journal={Duke Math. J.},
      volume={170},
      number={7},
       pages={1517\ndash 1583},
         url={https://doi.org/10.1215/00127094-2020-0074},
      review={\MR{4255058}},
}

\bib{Lindenstrauss2022}{unpublished}{
      author={Lindenstrauss, Elon},
      author={Wei, Daren},
       title={Time change rigidity for unipotent flows},
        date={2023},
        note={preprint},
}

\bib{margulis1971}{article}{
      author={Margulis, G.~A.},
       title={The action of unipotent groups in a lattice space},
        date={1971},
     journal={Mat. Sb. (N.S.)},
      volume={86(128)},
       pages={552\ndash 556},
      review={\MR{0291352}},
}

\bib{MargulisDiscrete1991}{book}{
      author={Margulis, G.~A.},
       title={Discrete subgroups of semisimple {L}ie groups},
      series={Ergebnisse der Mathematik und ihrer Grenzgebiete (3) [Results in
  Mathematics and Related Areas (3)]},
   publisher={Springer-Verlag, Berlin},
        date={1991},
      volume={17},
        ISBN={3-540-12179-X},
         url={https://doi.org/10.1007/978-3-642-51445-6},
      review={\MR{1090825}},
}

\bib{Morris05Ratner}{book}{
      author={Morris, Dave~Witte},
       title={Ratner's theorems on unipotent flows},
      series={Chicago Lectures in Mathematics},
   publisher={University of Chicago Press, Chicago, IL},
        date={2005},
        ISBN={0-226-53983-0; 0-226-53984-9},
      review={\MR{2158954}},
}

\bib{ornstein1982equivalence}{book}{
      author={Ornstein, Donald~S.},
      author={Rudolph, Daniel~J.},
      author={Weiss, Benjamin},
       title={Equivalence of measure preserving transformations},
      series={Mem. Amer. Math. Soc.},
        date={1982},
      volume={37},
      number={262},
         url={https://doi.org/10.1090/memo/0262},
      review={\MR{653094}},
}

\bib{ratner1978horocycle}{article}{
      author={Ratner, Marina},
       title={Horocycle flows are loosely {B}ernoulli},
        date={1978},
        ISSN={0021-2172},
     journal={Israel J. Math.},
      volume={31},
      number={2},
       pages={122\ndash 132},
         url={https://doi.org/10.1007/BF02760543},
      review={\MR{516248}},
}

\bib{ratner1979cartesian}{article}{
      author={Ratner, Marina},
       title={The {C}artesian square of the horocycle flow is not loosely
  {B}ernoulli},
        date={1979},
        ISSN={0021-2172},
     journal={Israel J. Math.},
      volume={34},
      number={1-2},
       pages={72\ndash 96 (1980)},
         url={https://doi.org/10.1007/BF02761825},
      review={\MR{571396}},
}

\bib{ratner1981some}{article}{
      author={Ratner, Marina},
       title={Some invariants of {K}akutani equivalence},
        date={1981},
        ISSN={0021-2172},
     journal={Israel J. Math.},
      volume={38},
      number={3},
       pages={231\ndash 240},
         url={https://doi.org/10.1007/BF02760808},
      review={\MR{605381}},
}

\bib{ratner1982rigidity}{article}{
      author={Ratner, Marina},
       title={Rigidity of horocycle flows},
        date={1982},
        ISSN={0003-486X},
     journal={Ann. of Math. (2)},
      volume={115},
      number={3},
       pages={597\ndash 614},
         url={https://doi.org/10.2307/2007014},
      review={\MR{657240}},
}

\bib{ratner1983horocycle}{article}{
      author={Ratner, Marina},
       title={Horocycle flows, joinings and rigidity of products},
        date={1983},
        ISSN={0003-486X},
     journal={Ann. of Math. (2)},
      volume={118},
      number={2},
       pages={277\ndash 313},
         url={https://doi.org/10.2307/2007030},
      review={\MR{717825}},
}

\bib{ratner1986rigidity}{article}{
      author={Ratner, Marina},
       title={Rigidity of time changes for horocycle flows},
        date={1986},
        ISSN={0001-5962},
     journal={Acta Math.},
      volume={156},
      number={1-2},
       pages={1\ndash 32},
         url={https://doi.org/10.1007/BF02399199},
      review={\MR{822329}},
}

\bib{ratner1990measure}{article}{
      author={Ratner, Marina},
       title={On measure rigidity of unipotent subgroups of semisimple groups},
        date={1990},
        ISSN={0001-5962},
     journal={Acta Math.},
      volume={165},
      number={3-4},
       pages={229\ndash 309},
         url={https://doi.org/10.1007/BF02391906},
      review={\MR{1075042}},
}

\bib{ratner1991measure}{article}{
      author={Ratner, Marina},
       title={On {R}aghunathan's measure conjecture},
        date={1991},
        ISSN={0003-486X},
     journal={Ann. of Math. (2)},
      volume={134},
      number={3},
       pages={545\ndash 607},
         url={https://doi.org/10.2307/2944357},
      review={\MR{1135878}},
}

\bib{ratner1991topological}{article}{
      author={Ratner, Marina},
       title={Raghunathan's topological conjecture and distributions of
  unipotent flows},
        date={1991},
        ISSN={0012-7094},
     journal={Duke Math. J.},
      volume={63},
      number={1},
       pages={235\ndash 280},
         url={https://doi.org/10.1215/S0012-7094-91-06311-8},
      review={\MR{1106945}},
}

\bib{ratner1995interactions}{inproceedings}{
      author={Ratner, Marina},
       title={Interactions between ergodic theory, {L}ie groups, and number
  theory},
        date={1995},
   booktitle={Proceedings of the {I}nternational {C}ongress of
  {M}athematicians, {V}ol. 1, 2 ({Z}\"{u}rich, 1994)},
   publisher={Birkh\"{a}user, Basel},
       pages={157\ndash 182},
      review={\MR{1403920}},
}

\bib{tang2020new}{article}{
      author={Tang, Siyuan},
       title={New time-changes of unipotent flows on quotients of {L}orentz
  groups},
        date={2022},
        ISSN={1930-5311},
     journal={J. Mod. Dyn.},
      volume={18},
       pages={13\ndash 67},
         url={https://doi.org/10.3934/jmd.2022002},
      review={\MR{4385813}},
}

\bib{witte1985rigidity}{article}{
      author={Witte, Dave},
       title={Rigidity of some translations on homogeneous spaces},
        date={1985},
        ISSN={0020-9910},
     journal={Invent. Math.},
      volume={81},
      number={1},
       pages={1\ndash 27},
         url={https://doi.org/10.1007/BF01388769},
      review={\MR{796188}},
}

\bib{witte1987zero}{article}{
      author={Witte, Dave},
       title={Zero-entropy affine maps on homogeneous spaces},
        date={1987},
        ISSN={0002-9327},
     journal={Amer. J. Math.},
      volume={109},
      number={5},
       pages={927\ndash 961},
         url={https://doi.org/10.2307/2374495},
      review={\MR{910358}},
}

\end{biblist}
\end{bibdiv}
\end{document}